%               This is an AmSTeX file
%====================================================================

\documentstyle{amsppt}

\TagsOnLeft \CenteredTagsOnSplits \NoRunningHeads
\topmatter
\title
Four drafts of the representation theory of the group of infinite
matrices over a finite field
\endtitle
\author A.~M.~Vershik, S.~V.~Kerov
\endauthor
\address
St.~Petersburg Department of
Steklov Institute of Mathematics, St.Petersburg, Russia
\endaddress
\email vershik\@pdmi.ras.ru \endemail
\endtopmatter

\document
\head Preface 2007. About these texts \endhead

The four texts presented below were written in 1997--2000, when
S.~V.~Kerov and me decided to return to the program
we had designed in the early 80s: to
generalize
the representation theory of the infinite symmetric group, which we
had started in the late 70s,
 to the groups of infinite matrices over finite fields.
Our first basic publication on this
subject \cite{VK} was unfortunately almost the last at that time ---
the only later publication was
my two introductory lectures published in the
Proceedings of the EMS Summer School \cite{V} (the
second lecture was devoted to this topic). The four
notes that are presented below have not been completely finished and are not
yet published. After Sergey's untimely death in July 2000, I did
not continue to work on this project, only gave several talks
in about 2001--2002 on our preliminary results. Now I have decided to keep
them in the original form, except for adding a few
additional remarks, corrections,
and references,\footnote{The list of references
given at the end of the fourth draft is common
for all four drafts and contains only
strictly necessary references.}
so this is a kind of
prepublication, whose goal is to stimulate and continue the
study of the subject which I consider extremely important,
especially nowadays, because of the increasing interest of
mathematicians to various $q$-analogs of classical objects.

It makes sense to read these drafts after the paper \cite{VK}, but we
will not reproduce it here. Nevertheless, we recall some important
facts that will help to join these drafts together.

It is well known that the main object that can be considered
as a $q$-analog of the symmetric group ${\frak S}_n$ is the group
$GL(q,n)$ of matrices over the finite field $\Bbb F_q$ (with $q$ a power of a prime).
   But it turned out that for the infinite symmetric group ${\frak S}_{\infty}$,
   the right $q$-analog is not the group $GL(q,\infty)$, which is the
   inductive limit of the groups $GL(q,n)$ in $n$, but another group,
   which we denoted by $GLB(\Bbb F_q)$; it is neither
   inductive nor projective limit of finite groups, but a
   so-called IP-limit (inductive-projective limit), see below.
   This is the locally compact group of all infinite matrices over the finite
   field that are finite (= contain finitely many nonzero elements) below
   the main diagonal.
   Thus the Borel subgroup of this group is compact. At the same time,
   one version of the group algebra of the
   group $GLB$ (an analog of the Bruhat--Schwartz algebra of the classical $p$-adic
   groups) is the inductive
   limit of the group algebras ${\Bbb C}(GL(q,n))$,
   so that we can apply
   the methods of inductive theory. This inductive limit of algebras
   does not correspond to any inductive limit of the corresponding groups; at
   the group level, this means that we use a more complicated
   construction, which is described here. The reason why the group $GLB$
   is the right infinite $q$-analog of ${\frak S}_{\infty}$ is related to
   the fundamental fact of the classical representation theory
   of the groups $GL(q,n)$:
   if we want the
   restrictions of irreducible representations of
   $GL(q,n+1)$ to $GL(q,n)$ to be multiplicity-free,
   we must use
   the {\it parabolic embedding} of $GL(q,n)$
   into $GL(q,n+1)$
   (see \cite{Green, Fad, Zel}).
    More exactly, the natural procedure for extending a
   representation of $GL(q,n)$ to a representation of
   $GL(q,n+1)$ is not the usual induction from $GL(q,n)$
   as a subgroup of $GL(q,n+1)$, but the induction from a maximal parabolic
   subgroup of $GL(q,n+1)$ that contain $GL(q,n)$.
   Thus we need to respect this rule for embedding
   and generalize it to the infinite case.\footnote{As mentioned in
   \cite{VK}, the notion of the group $GLB$ was suggested by the authors and
   A.~Zelevinsky in 1983.}

   In the first draft presented below
   (``Interpolation between inductive and projective limits of finite
groups with applications to linear groups over finite fields'')
we consider
   a general scheme called the IP-limit of groups,
   the main example of which is the group $GLB$ as the IP-limit of $GL(q,n)$.
   It is very interesting to find other nontrivial examples
   of this scheme.
It is worth adding that the idea of the mixed (indo-projective, or IP)
   limit was used in several categories (convex sets, linear spaces,
   etc.). The first author used it for
   studying the notion of
   the tower of measures (see \cite{V94}). Here we
   consider the pure group-theoretical aspect of this notion.

The next note, ``The characters of the group of almost triangular
matrices over a finite field,'' is devoted to the main topic, the theory of
characters of the groups under consideration. First of all,
at the first stage it is better
to restrict ourselves to the so-called principal series of
representations and unipotent characters.
It is worth mentioning that
\cite{VK} contains a mistake, which was found after recent
discussions with E.~Goryachko. Namely, the Bratteli diagram for the
group $GLB$ is not the disjoint union of copies of the Young graph,
as we
claimed in \cite{VK}; fortunately, this claim does not affect any other
assertion in \cite{VK}. The diagram in question is indeed the
disjoint union of
countably many components, but each of these components is the direct
product of copies of the Young graph (depending on the field). Thus the
principal series is only a part of the main component of the
Bratteli diagram; nevertheless,
the branching rule and the structure of characters
of all other series (caspidal etc.) are almost the same as
for the principal series.
That is why we switch on from the group $GLB$
to the subgroup $GLN(q)$, which is the IP-limit with the ``unipotent''
projections; the difference between these groups
is merely that the diagonal elements of matrices in $GLN$
are eventually 1's
and the maximal unipotent subgroup
consists of upper triangular
matrices with 1's on the diagonal.

We formulate a series of conjectures about the unipotent
characters of the group $GLN$, and the main fact is (see also the
third draft) that every indecomposable character determines on the Borel
(A.~Borel!) subgroup {\it a true Borel} (E.~Borel!) probability
measure, and these measures are mutually singular for different characters.

In the next draft, ``A Law of Large Numbers for the characters of
$GL_n(k)$ over a finite field $k$,'' which partially reproduces the
previous one, we try to suggest a program of how to connect the theory of
characters with the structure of Jordan forms of typical
 matrices (in the sense of the above-mentioned measures).
 This note concerns a deep
 and very intriguing link
 between two appearances of the Thoma parameters
 (see \cite{Thoma}): in the first
 case they arise in the representation theory of
 the infinite symmetric group
 and the Hecke algebras as the frequences of rows and columns of
 Young diagrams
 corresponding to representations; and in the second case they arise as
 the normalized lengths of Jordan blocks of matrices. This is a
 manifestation
 of some kind of duality,  which also appeared in our works on models
 of factor representations of the infinite symmetric group.
In our theory the characters are realized as the traces of well-defined
type $II_{\infty}$ representations.
Note that this is perhaps the first occasion when an infinite trace on
the group algebra of a locally finite group can be interpreted not as a
function on the group
and not as a degenerate linear functional on the group algebra, but as
a singular
(with respect to the Haar measure)  measure on the group; moreover,
different characters generate mutually singular measures, and
for most characters the measures are supported by
the whole group.

The last draft is devoted to the
problem concerning realizations of representations, which is
 most important for applications.
We give the simplest example of an implicit construction of
principal series representations of the group $GLB$,
which is based on the geometry of the
infinite-dimensional Grassmanian; from the point of view of Schubert
cells, this construction corresponds
to infinite Young diagrams with two rows.
One of the puzzles here is how a type $III$ factor representation
comes to the picture. It seems to me that construction
of factor representations and some irreducible representations
of GLB-type groups is
the mainstream of the forthcoming research in this area, which
will be continued.

I am deeply grateful to N.~V.~Tsilevich for her valuable help
in preparing the manuscript for publication.

\bigskip
\hfill A.~Vershik

\newpage

\topmatter
\title
I. Interpolating between inductive and projective limits of finite
groups with applications to linear groups over finite fields
\endtitle

\endtopmatter

\subhead 1. IP-families of finite groups \endsubhead In this section
we introduce IP-families, generalizing both inductive and projective
families of finite groups.

\smallskip\noindent

{\bf Definition.} Consider a system
$$
\{1\} = P_1 \subset G_1 \leftarrow
P_2 \subset G_2 \leftarrow
P_3 \subset G_3 \leftarrow
\ldots
$$
of finite groups and group homomorphisms, and assume that all
maps $\pi:\,P_{m+1}\to G_m$ are epimorphisms. We refer to such
systems as {\it IP-families}.

\smallskip\noindent
{\bf Example 1.}
Every inductive family
$$
G_1 \subset G_2 \subset G_3 \subset \ldots
$$
determines a trivial IP-family with $P_{m+1}=G_m$.

\smallskip\noindent
{\bf Example 2.}
Every projective family
$$
G_1 \leftarrow G_2 \leftarrow G_3 \leftarrow \ldots
$$
determines a trivial IP-family with $P_m=G_m$.

The following two examples are less trivial and motivate our
definition of IP-families.

\smallskip\noindent
{\bf Example 3.}
Let $k=\Bbb F_q$ be the finite field with $q$ elements and
$G_m=GL_m(k)$ the general linear group over $k$. We take for
$P_{m+1}\subset G_{m+1}$ the subgroup of block matrices of the
form
$$
g = \pmatrix A & b \\ 0 & a \endpmatrix
$$
where $A\in G_m$, $a\in G_1\cong k^*$, and $b\in k^m$. The
epimorphism $\pi:P_{m+1}\to G_m$ maps $g\in P_{m+1}$ to the
matrix $A$. This is our basic example of an IP-family.

\smallskip\noindent
{\bf Example 4.}
The groups $G_m=GL_m(k)$ are the same as in Example 3, but we
take for $P_{m+1}$ the affine subgroup of matrices of the form
$$
g = \pmatrix A & b \\ 0 & 1 \endpmatrix.
$$
The epimorphisms $\pi:P_{m+1}\to G_m$ are defined as in Example
3.

\smallskip\noindent
{\bf Example 5.}
Here is a baby example intermediate between the trivial examples 1, 2
and the main examples 3, 4. Let $H$ be a finite group, and
denote by $G_m$ the wreath product of the symmetric group $\frak
S_m$ with the coefficient group $H$. One can think of $G_m$ as
the group of $m\times m$ permutation matrices with nontrivial
values in the group $H$. We take for $P_{m+1}$ the subgroup of
matrices in $G_{m+1}$ with a nontrivial element on the crossing
of the last row and the last column (so that the other elements
of the last row and column vanish). The value of the
epimorphism $\pi$ at a matrix $g\in P_{m+1}$ is, by definition,
the submatrix of $g$ on the crossing of the first $m$ rows and
columns.

\subhead 2. Associated limiting groups \endsubhead
In this section we define the limiting group of an IP-family, and
its basic compact subgroup.

We start with a definition of subgroups $G(m,n)\subset G_n$, for
all $1\le m\le n$. To this end we set $G(m,m)=G_m$, and we define
$G(m,n+1)$ by induction as the preimage
$\pi^{-1}\big(G(m,n)\big)$ of the group $G(m,n)\subset G_n$ in
$P_{n+1}\subset G_{n+1}$. In particular, $G(m,m+1)=P_{m+1}$.

The restrictions of the homomorphisms $\pi$ to the subgroups
$G(m,n)\subset P_n$ determine epimorphisms
$\pi:G(m,n+1)\to G(m,n)$, so that a projective family
$$
G(m,m) \leftarrow
G(m,m+1) \leftarrow
G(m,m+2) \leftarrow \ldots
$$
arises. We denote by $G(m,\infty)=\varprojlim_n G(m,n)$ the
corresponding profinite group. By definition, an element $g\in
G(m,\infty)$ is a sequence $g=\{g_n\}_{n=m}^\infty$ of elements
$g_n\in G(m,n)$ such that $\pi(g_{n+1})=g_n$ for all $n\ge m$.

Given an element $g=\{g_n\}_{n=m}^\infty\in G(m,\infty)$, the
truncated sequence $\widetilde{g}=\{g_n\}_{n=m+1}^\infty$
determines an element $\widetilde{g}$ of the group
$G(m+1,\infty)$. It follows that the groups $G(m,\infty)$ form an
inductive family of compact groups
$$
G(1,\infty) \subset
G(2,\infty) \subset
G(3,\infty) \subset \ldots.
$$

\smallskip\noindent
{\bf Definition.}
We denote by
$$
G = \varinjlim_m \varprojlim_n G(m,n)
$$
the inductive limit of the groups $G(m,\infty)$, and we consider
$G$ as the {\it limiting group} of an IP-family. In the inductive
limit topology the group $G$ is locally compact and totally
disconnected. The compact subgroup $B=G(1,\infty)\subset G$ will
be referred to as the {\it basic profinite subgroup} of the group
$G$.

We denote by $\mu$ the Haar measure on the group $G$ normalized
by the condition $\mu(B)=1$.

\subhead 3. The Bruhat--Schwartz group algebra \endsubhead
Let $\Cal B$ denote the linear space of locally constant finitely
supported functions on $G$ with complex values. The space $\Cal
B$ is closed under the convolution
$$
(f_1*f_2)\,(g) = \int_G f_1(gh^{-1})\,f_2(h)\, d\mu(h), \qquad
f_1, f_2 \in \Cal B,
$$
and under the standard involution
$$
f^\#(g) = \overline{f(g^{-1})}, \qquad f \in \Cal B.
$$
The functions in $\Cal B$ separate the points of the group $G$.
We use the $*$-algebra $\Cal B$ as the basic group algebra of the
group $G$.

The main advantage of this choice of $\Cal B$ is that this
algebra is locally semisimple, hence it can be studied by the
powerful combinatorial techniques of Bratteli diagrams (see
\cite{LSS}).

\define\Cyl{\operatorname{Cyl}}

Let $\Cal B_m$ be the $*$-subalgebra of those functions in $\Cal
B$ that are supported by the subgroup $G(m,\infty)\subset G$ and does not
depend but on the image $g_m\in G_m$ of an element
$g=\{g_n\}_{n=m}^\infty\in G(m,\infty)$. The algebra $\Cal B_m$
is clearly $*$-isomorphic to the group algebra $\Bbb C[G_m]$. Let
$\widetilde{b}_g\in \Cal B_m$ denote the characteristic function
of the cylinder subset $\Cyl(g)\subset G(m,\infty)$, where
$$
\Cyl(g) = \{g=\{g_n\}_{n=m}^\infty\in G(m,\infty):\; g_m=g\},
\qquad g \in G_m.
$$
We also set $b_g=\widetilde{b}_g/|N_m|$, where $N_m$ is the
(finite) kernel of the epimorphism $\pi:\,P_{m+1}\to G_m$. Then
$$
b_g*b_h = b_{gh}, \qquad (b_g)^\# = b_{g^{-1}},
\qquad g,h \in G_m.
$$
If $E_m$ denotes the identity element of the group $G_m$, then
$b_{E_m}$ is the identity in the algebra $\Cal B_m$.

\proclaim{Proposition}
The algebra $\Cal B$ is locally semisimple. More precisely, $\Cal
B = \varinjlim \Cal B_m$ is the limit of the inductive family
$$
\Cal B_1 \subset
\Cal B_2 \subset
\Cal B_3 \subset \ldots
\tag$*$
$$
of semisimple finite dimensional subalgebras $\Cal B_m \subset
\Cal B$.
\endproclaim

\demo{Proof}
In terms of the basis elements $b_g$ the algebra inclusion
$\Cal B_m\subset\Cal B_{m+1}$ is determined by the formula
$$
b_g \mapsto b_{i(g)} =
{1\over |N_m|} \sum_{h\in N_m} b_{gh}.
$$
Here we use the notation
$$
i(g) = {1\over |N_m|} \sum_{h\in N_m} gh
$$
for the corresponding homomorphism $i:\,\Bbb C[G_m]\to\Bbb
C[G_{m+1}]$. One can readily check that this is indeed a
homomorphism of $*$-algebras, and that $\Cal B$ is the union of
the subalgebras $\Cal B_m$, $m\ge1$. Note that the map $i$ does
not respect the units in the algebras $\Cal B_m$. As a result,
there is no unit in the algebra $\Cal B$.
\qed\enddemo

\subhead 4. The branching of characters \endsubhead
In this section we describe the Bratteli diagram of the inductive
family $\{\Cal B_m\}$ approximating the algebra $\Cal B$.

Let $\widehat{G}_m$ denote the finite set of equivalence classes
of irreducible complex representations of the group $G_m$. We
take $\Gamma_m=\widehat{G}_m$ for the $m$th level of the
branching diagram $\Gamma=\bigcup_{m=0}^\infty\Gamma_m$.

\define\Ind{\operatorname{Ind}}

Given an irreducible representation $\varphi\in\Gamma_m$, we
denote by $\widetilde{\varphi}$ the irreducible representation of
the group $P_{m+1}$ in the same Hilbert space lifted from
$\varphi$ via the epimorphism $\pi:\,P_{m+1}\to G_m$. Let
$$
\Ind\widetilde{\varphi} =
\bigoplus_\psi \varkappa(\varphi,\psi)\;\psi
$$
be the decomposition into irreducibles of the representation
$\Ind\widetilde{\varphi}$ induced to $G_{m+1}$ by the
representation $\widetilde{\varphi}$ of $P_{m+1}$. We denote by
$\varkappa(\varphi,\psi)\in\Bbb N$ the multiplicity of $\psi$ in
the decomposition.

\smallskip\noindent
{\bf Definition.}
Let $\Gamma=\bigcup_{m=0}^\infty\Gamma_m$ be the vertex set of a
branching diagram, and define the multiplicities of transitions
as $\varkappa(\varphi,\psi)$, where $\varphi\in\Gamma_m$ and
$\psi\in\Gamma_{m+1}$, for all $m\ge0$. This is the Bratteli
diagram of the inductive system ($*$). We shall refer to the
branching graph $(\Gamma,\varkappa)$ as the {\it character
branching} of the IP-family.

The character branching describes the restrictions of characters
of irreducible representations of the algebras $\Cal B_m$.
Choose a representation $\psi\in\Gamma_{m+1}$ and consider the
restriction $\chi^\psi(i(g))$ of its character
$\chi^\psi(g)=\operatorname{Tr}\psi(g)$ to the elements $g\in G_m$.
Then
$$
{1\over|N_m|} \sum_{h\in N_m} \chi^\psi(gh) =
\sum_{\varphi\in\Gamma_m} \varkappa(\varphi,\psi) \chi^\varphi(g).
\tag$**$
$$

\smallskip
\noindent{\bf Remark.} {The branching diagram $(\Gamma,\varkappa)$ determines
the inductive family of algebras $\{\Cal B_m\}$ up to the dimensions
of its irreducible representations. The reason is that we disregard
the zero components arising in the decompositions of restricted
representations. The limiting algebra $\Cal B$ is determined by
$(\Gamma,\varkappa)$ up to stable isomorphism.}
\medskip

Let $\Delta_m\subset\Gamma_m$ denote the set of irreducible
representations $\varphi$ of $G_m$ which contain a nonzero
$G(1,m)$-invariant vector. Representations in
$\Delta=\bigcup_{m=0}^\infty\Delta_m$ will be referred to as {\it
unipotent} irreducible representations of the group $G_m$. Note
that every nonzero irreducible component of the restriction
$\psi(i(\cdot))$ of a unipotent representation
$\psi\in\Delta_{m+1}$ is itself a unipotent representation in
$\Delta_m$, so that $(\Delta,\varkappa)$ is a branching graph,
too. In fact, let $\sigma_m$ denote the representation of $G_m$
induced by the identity representation of the subgroup $G(1,m)$.
(In our Example 3 this is the {\it flag} representation).
The unipotent representations
can be alternatively defined as the
irreducible components of $\sigma_m$. The claim now follows from
the observation that $\sigma_m$, when extended trivially to
$P_{m+1}$ and induced to $G_{m+1}\supset P_{m+1}$, is equivalent
to $\sigma_{m+1}$. The limiting algebra corresponding to the
branching graph $(\Delta,\varkappa)$ is the Hecke algebra of the
pair $G\supset B$.

\subhead 5. The branching of conjugacy classes \endsubhead
There is another branching graph associated with every IP-family,
and we proceed with the definition of this graph.

Let us say that elements $g_1,g_2\in G(1,m)\subset G_m$ are
{\it equivalent}, or $G$-{\it conjugate}, if $g_2=h^{-1}g_1h$ for
some $h\in G_m$. In this case we write $g_1\sim g_2$.

Given a probability measure $M$ on the group $B=\varprojlim
G(1,m)$, we denote by $M_m(g)$ the cylinder probabilities
$M_m(g) = M\{\Cyl(g)\}$, $g\in G(1,m)$. The coherence condition
reads as
$$
M_m(g) = \sum_{h\in N_m} M_{m+1}(gh) \qquad
\text{ for all } g\in G(1,m).
$$

\medskip\noindent
{\bf Definition.}
We say that a probability measure $M$ on the group $B$ is
$G$-{\it central} if its cylinder probabilities $M_m(g)$
depend only on the $G_m$-conjugacy class of $g\in G(1,m)$, for all
$m\ge1$. A $G$-central measure $M$ is {\it ergodic} if it
cannot be represented as a nontrivial mixture of other $G$-central
measures.

\medskip\noindent
{\bf Example.}
The normalized Haar measure of the group $B$ is an example of a
$G$-central measure, for every IP-family.

\medskip\noindent
{\bf Problem.}
Describe all ergodic $G$-central measures on the basic profinite
subgroup $B$ of a given IP-family.
\medskip

Assume that a character $\chi$ of the group $G$ (or, more
precisely, of the Bruhat--Schwartz group algebra $\Cal B$) is
nonnegative on the subgroup $B$. Then the identity
$$
M_m(g) = \chi(\widetilde{b}_g), \qquad g \in G(1,m),
$$
defines a probability distribution on the group $G(1,m)$, for all
$m\ge1$. By the formula ($**$), the family of distributions
$\{M_m\}$ is coherent, hence determines a $G$-central measure on
$B$.

\smallskip\noindent
{\bf Example.}
Consider the IP-family of Example 5, Section 1 assuming for
simplicity that the coefficient group $H$ is Abelian. The group
$G$ in this example is the wreath product of the group $\frak
S_\infty$ of all finite permutations of the set $\Bbb N$ with
coefficients in $H$. The subgroup $B$ is the subgroup of diagonal
matrices, hence it is isomorphic to the infinite direct product
$B\cong H^\infty\equiv\prod_1^\infty H$.

In this example the subgroup $B$ is normal in $G$. It follows
that the only unipotent representations are those which are
trivial on $B$. All of them induce the Haar measure on the group
$B$.

By the celebrated de Finetti Theorem, the class of $G$-central
measures on $B$ coincides with that of product measures
$M=\prod_1^\infty M_0$, for some probability distribution $M_0$
on $H$. For instance, take for $M_0$ the Haar measure on a
subgroup $H_0\subset H$. Then the corresponding $G$-central
measure is, in a sense, positive definite, though it cannot be
represented as a mixture of {\it characteristic measures} arising
from the characters of $G$, nonnegative on $B$.

\newpage

\topmatter
\title
II. The characters of the group of almost triangular matrices over a
finite field
\endtitle
\endtopmatter

This paper is a sequel of \cite{VK}. We continue the study of
the asymptotic character theory of the general linear group
$GL(n,q)$ over a finite field $k=\Bbb F_q$ with $q$ elements.
Here we focus on unipotent characters, and our goal is basically
to describe the asymptotics of character values at unipotent
conjugacy classes.

Our approach to the asymptotic character theory of the groups
$GL(n,q)$ is governed by the character theory of the
infinite-dimensional limiting group $GLU(q)$. A character of this group is
entirely determined by a probability measure on the group $U$ of
upper unitriangular matrices, and we study the statistical
properties of such measures.

The group $GLU(q)$ is similar to the group $GLB(q)$ introduced in
\cite{VK}. Both $GLB(q)$ and $GLU(q)$ have much more characters than
the group $GL(\infty,q)$ studied by E.~Thoma \cite{Th}
and H.~Skudlarek \cite{Sk}.
It is this
phenomenon that motivates our way of embedding the group algebras
$\Bbb C[GL(n,q)]$. Another motivation comes from the fact that the
branching of unipotent characters with respect to our embedding is
governed by the Young lattice, hence coincides with the branching of
irreducible characters of the symmetric groups.

We refer to \cite{LSS} for the character theory of locally
semisimple algebras, to \cite{Mac} for the definitions of symmetric
functions, and to \cite{Zel} for the representation theory of the
groups $GL(n,q)$.

\subhead 0. Unipotent characters of $GL(n,q)$
\endsubhead
Let $B$ be the group of upper triangular matrices, and denote by
$\Cal F=GL(n,q)/B$ the flag space. Irreducible characters
$\chi^\lambda$ contained in the module $\Bbb C[\Cal F]$ are called
{\it unipotent}. They are labeled by Young diagrams with $n$ boxes.
The degree of $\chi^\lambda$ is given by the following $q$-version
of the {\it hook formula}:
$$
\chi^\lambda(E) = q^{n(q)}
\prod_{j=1}^n (q^j-1) \prod_{b\in\lambda} (q^{h(b)}-1)^{-1},
$$
where $n(\lambda)=\sum_{k\ge1}(k-1)\lambda_k$. Note that
$\chi^{(n)}$ is the identity character, and $\chi^{(1^n)}$ is
known as the Steinberg character.

Let $P_\mu$ be a parabolic subgroup of type $\mu\in \Bbb Y_n$. We
denote by $\psi^\mu$ the character of the natural representation
in $\Bbb C[GL(n,q)/P_\mu]$. The decomposition of $\psi^\mu$ into
irreducible characters has the form
$$
\psi^\mu = \sum_\lambda K_{\lambda\mu}\, \chi^\lambda,
$$
where $K_{\lambda\mu}$ are known as the {\it Kostka numbers}. Note
that the coefficients $K_{\lambda\mu}$ do not depend on $q$.

\subhead 1. The group $GLU(q)$ \endsubhead
Let $k=\Bbb F_q$ be the finite field with $q$ elements (the field
is known to be unique up to isomorphism). An infinite
matrix $g=(g_{ij})_{i,j=1}^\infty$ over the field $k$ is called
{\it almost unitriangular} if

\roster
\item there are only finitely many nonzero elements below the
main diagonal;

\item all but finitely many diagonal elements are equal to one.
\endroster

The group $GLU(q)$ of all invertible almost unitriangular
matrices is a locally compact topological group. The group $U$ of
upper unitriangular matrices is a compact open subgroup in
$GLU(q)$. The topology of $U$ is that of the product of the finite
set $k$ over all matrix elements above the diagonal. These matrix
elements are independent and uniformly distributed in $k$ with
respect to the normalized Haar measure $M$ of $U$.

\subhead 2. The Schwartz--Bruhat group algebra \endsubhead
The group algebra $\Cal B$ of locally constant compactly
supported functions on the group $GLU(q)$ can be naturally
represented as the limit of the inductive family
$$
\Bbb C[GL_1(k)] \subset \Bbb C[GL_2(k)] \subset \ldots \subset
\Bbb C[GL_n(k)] \subset \ldots
\tag 1
$$
of the complex group algebras of the groups $GL_n(k)$. The group
$GL_\infty(k)=\varinjlim GL_n(k)$ is a countable dense subgroup
in $GLU(q)$. Note that the family (1) is {\it not} generated by
the natural embeddings of the groups $GL_n(k)$.

\subhead 3. Characteristic measures \endsubhead
Every positive definite class function $\chi:GLU(q)\to\Bbb C$ is
called a {\it character} of $GLU(q)$. A character $\chi$ is {\it
indecomposable} if every decomposition $\chi=\chi_1+\chi_2$ in a
sum of two characters is trivial in the way that both terms are
constant multiples of $\chi$.

We find all indecomposable characters of the group $GLU(q)$. In
particular, there is a family of {\it unipotent} indecomposable
characters $\chi^{(\alpha;\beta)}$ labelled by pairs
$(\alpha;\beta)$ of nonnegative weakly decreasing sequences
$\alpha = (\alpha_1,\alpha_2,\ldots)$, $\beta =
(\beta_1,\beta_2,\ldots)$ such that
$\sum_{k\ge1}\alpha_k+\sum_{k\ge1}\beta_k\le1$. The space
$\Delta$ of such sequences is known as the {\it Thoma simplex}. The
parameter space $\Delta$ of indecomposable unipotent characters
is exactly the same as that of indecomposable characters of
the infinite symmetric group. Remark that the group
$GL_\infty(k)=\varinjlim GL_n(k)$ studied by E.~Thoma has only
countably many indecomposable characters.

\subhead 4. Central distributions of unitriangular matrices
\endsubhead
Let $U\subset GLU(q)$ be the subgroup of upper unitriangular
matrices in $GLU(q)$. Denote by $g^{(n)}$ the submatrix on the
crossing of the first $n$ rows and columns of a matrix $g\in U$,
and write $\rho^{(n)}$ for the Young diagram describing the
Jordan form of $g^{(n)}$. Denote by $U_h$ the cylinder set of all
matrices $g\in U$ with $g^{(n)}=h$, and by $\rho=\rho^{(n)}$ the
Jordan partition of $h$.

A probability measure $m$ on $U$ is called {\it central} if its
cylinder probabilities $m(U_h)$ depend on the Jordan partition
$\rho$ of the matrix $h$ only. We denote the cylinder
probabilities of a central measure as $m_\rho=m(U_h)$.

We say that a set $S$ of matrices in $U$ is a {\it tail set} if
$S$ contains, with every matrix $g\in S$, all the matrices $g'\in
U$ that differ from $g$ in only finitely many matrix elements. A
central measure $m$ is called {\it ergodic} if every measurable
tail set has probability $0$ or $1$.

Let $Q_\rho(\alpha;\beta;t)$ denote the dual Hall--Littlewood
(super) symmetric function with parameter $t$.

\proclaim{Theorem}
For every point $(\alpha;\beta)\in\Delta$ of the Thoma simplex
there exists an ergodic central measure $m^{(\alpha;\beta)}$ with
cylinder probabilities
$$
m^{(\alpha;\beta)}_\rho =
{Q_\rho(\alpha;\beta;t) \over t^{n(\rho)} (1-t)^n},
\tag 2
$$
where $t=1/q$, and $n(\rho)=\sum_{k\ge1} (k-1)\rho_k$ for a
partition $\rho=(\rho_1,\rho_2,\ldots)\vdash n$.
\endproclaim

\proclaim{Conjecture 1}
The above formula provides all ergodic central measures on $U$.
\endproclaim

Given a matrix $g\in U$, the sequence
$t=(\rho^{(1)},\rho^{(2)},\ldots)$ is an infinite Young tableau.
Set $r_k(\rho)$ for the length of the $k$th row of a Young
diagram $\rho$, and $c_k(\rho)$ for the length of the $k$th
column of $\rho$.

\proclaim{Conjecture 2}
The limits
$$
\lim_{n\to\infty} {r_k(\rho^{(n)}) \over n} = \alpha_k, \qquad
\lim_{n\to\infty} {c_k(\rho^{(n)}) \over n} = \beta_k
\tag 3
$$
exist for almost all matrices $g\in U$ with respect to the
measure $m^{(\alpha;\beta)}$.
\endproclaim

\subhead 5. Characteristic measures of unipotent characters
\endsubhead
Every unipotent character $\chi^{(\alpha;\beta)}$ of the group
$GLU(q)$ determines a central probability distribution
$M^{(\alpha;\beta)}$ on the group $U$ referred to as its {\it
characteristic measure}. We show that
$$
M^{(\alpha;\beta)}(U_h) =
r_\rho(\alpha;\beta;t) \equiv t^{-n(\rho)}
\sum_{\lambda\vdash n} K_{\lambda\rho}(t)\;
s_\lambda(\alpha;\beta),
\tag 4
$$
where $t=1/q$, $K_{\lambda\rho}(t)$ is the Kostka--Foulkes matrix,
and $s_\lambda$ is the Schur (super) symmetric function.

In order to identify the characteristic measure
$M^{(\alpha;\beta)}$ of a unipotent indecomposable character
$\chi^{(\alpha;\beta)}$ as one of the measures
$m^{(\alpha;\beta)}$ we need more notation. Given $0<t<1$ and a
positive number $a$, we denote by $a^{(t)}$ the geometric
sequence
$$
a^{(t)} = (1-t)a,(1-t)ta,(1-t)t^2a,\ldots
$$
If $\alpha=(\alpha_1,\alpha_2,\,\ldots)$ is a sequence of
positive numbers, we denote by $\tilde\alpha$ the sequence
obtained from $\alpha$ by rearranging the elements of the
sequences $\alpha^{(t)}_1,\alpha^{(t)}_2,\,\ldots$ in decreasing
order.

\proclaim{Theorem}
The characteristic measure $M^{(\alpha;\beta)}$ of a unipotent
character $\chi^{(\alpha;\beta)}$ coincides with the central
measure $m^{(\tilde\alpha;\beta)}$ corresponding to the point
$(\tilde\alpha;\beta)\in\Delta$.
\endproclaim

\subhead 6. Extending unipotent characters from the subgroup $U$
\endsubhead
Every indecomposable unipotent character $\chi^{(\alpha;\beta)}$
of the group $GLU(q)$ is determined, as a functional on the group
algebra $\Cal B$, by a signed Radon measure $\tilde
M^{(\alpha;\beta)}$ on $GLU(q)$. Given a matrix $h\in GL_m(k)$,
consider the cylinder set $C_h\subset GLU(q)$ of all matrices $g\in
GLU(q)$ such that $g^{(m)}=h$ and $g_{ij}=0$ if $i>m$ and $i>j$.
The cylinder measures $\tilde M^{(\alpha;\beta)}(C_h)$ depend on
the conjugacy class of $h$ in $GL_n(k)$ only. Recall that the
conjugacy classes are labelled by the families $\varphi:\Cal
F(q)\to\Bbb{Y}$ of Young diagrams such that
$$
||\varphi|| \equiv \sum_{f\in\Cal F(q)} |\varphi(f)|\,d_f = m.
\tag 5
$$
Here $\Cal F(q)$ is the set of irreducible polynomials in $\Bbb
F_q[t]$ with the exception of the polynomial $t$, and $d_f$
denotes the degree of $f$.

We show that the cylinder measures $\tilde
M^{(\alpha;\beta)}(C_h)$ can be restored from the cylinder
probabilities of the corresponding characteristic measure
$M^{(\alpha;\beta)}$.

\proclaim{Theorem}
The cylinder measures $\tilde M^{(\alpha;\beta)}(C_h)$ can be
restored from the cylinder probabilities {\rm(4)} of the corresponding
characteristic measure $M^{(\alpha;\beta)}$ by the formula
$$
\tilde M^{(\alpha;\beta)}(C_h) =
\prod_{f\in\Cal F(q)} r_{\varphi(f)}
(\alpha^{d_f};(-1)^{d_f+1}\beta^{d_f};t^{d_f}).
\tag 6
$$
Here $\varphi$ is the family of Young diagrams describing the
conjugacy class of the matrix $h$, and we define the $d$th power
of a sequence $\alpha=(\alpha_1,\alpha_2,\ldots)$ as
$\alpha^d=(\alpha_1^d,\alpha_2^d,\ldots)$.
\endproclaim

\newpage

\topmatter
\title
III. A Law of Large Numbers for the characters of $GL_n(k)$
over a finite field $k$
\endtitle

\endtopmatter

\head Introduction \endhead

Let $k=\Bbb F_q$ denote the finite field with $q$ elements, and
consider the group $GL_\infty(k)$ of infinite nonsingular
matrices over $k$ of the form $g=E+f$, where $E$ is the identity
matrix and $f$ is a matrix with finitely many nonzero elements.
The characters of the group $GL_\infty(k)$ were studied and
classified by E.~Thoma and his successors (see \cite{Th, Sk}).
Since the group $GL_\infty(k)$ is the inductive limit
of the finite linear groups $GL_n(k)$, $n\to\infty$, characters
of $GL_\infty(k)$ can be regarded as ``coherent'' families of
characters of irreducible representations of the groups
$GL_n(k)$.

In a recent paper \cite{VK}
the authors have introduced a new
locally compact group, $GLB(k)$, made of all nonsingular infinite
matrices over the field $k$ with only finitely many nonzero
elements below the main diagonal. In a sense, this group can also
be regarded as a limit of the finite linear groups $GL_n(k)$, and its
characters are represented by coherent (in a sense to be made
precise below) families of irreducible characters of the latter
groups. One of the advantages of the group $GLB(k)$, compared
with $GL_\infty(k)$, is that it has much bigger amount of
characters, and that the formal substitution $q=1$ admits a
natural interpretation in terms of the character theory of the
infinite symmetric group $\frak S_\infty$.

The present paper is a sequel of \cite{VK}. Here we restrict
ourselves to unipotent characters of the groups $GL_n(k)$, and
the limits thereof. Basically, our claims are as follows:

\roster
\item Unipotent indecomposable characters
$\chi^{(\alpha;\beta)}$ of the group $GLB(k)$ (labelled by
elements $(\alpha;\beta)$ of the Thoma simplex, see \cite{VK})
can be treated as (signed) measures on this group.

\item Every such measure is uniquely determined by its
restriction $M^{(\alpha;\beta)}$ to the compact subgroup
$U(k)\subset GLB(k)$ of upper triangular matrices with the unit
diagonal. The restrictions, referred to as {\it characteristic
measures}, are in fact probability measures on $U(k)$.

\item Characteristic measures corresponding to different points
of the Thoma simplex are mutually singular.

\item We state a conjectural Law of Large Numbers for the
characteristic measures. In order to describe it in detail, we
need some notation. Let $M^{(\alpha;\beta)}$ be the
characteristic measure corresponding to Thoma parameters
$(\alpha;\beta)$. Let $u(n)$ denote the submatrix of a matrix
$u\in U(k)$ on the crossing of the first $n$ rows and columns,
and let $\lambda(n)$ be the partition of $n$ corresponding to the
Jordan block structure of $u(n)$. A matrix $u\in U(k)$ is said to
be {\it regular} if the following limits (called the row and column
{\it frequencies}) exist for all $k=1,2,\ldots$:
$$
\lim_{n\to\infty} {\lambda(n)_k \over n} = \widetilde{\alpha}_k,
\tag 1.1
$$
$$
\lim_{n\to\infty} {\lambda'(n)_k \over n} = \widetilde{\beta}_k.
\tag 1.2
$$
\endroster

\proclaim{Conjecture}
Almost all matrices $u\in U(k)$ are regular with respect to every
characteristic measure $M$. Set $t=1/q$ and assume that
$M=M^{(\alpha;\beta)}$. Then the list
$$
\aligned
\{&(1-t)\alpha_1,(1-t)t\alpha_1,(1-t)t^2\alpha_1, \ldots,\\
&(1-t)\alpha_2,(1-t)t\alpha_2,(1-t)t^2\alpha_2, \ldots,\\
&(1-t)\alpha_3,(1-t)t\alpha_3,(1-t)t^2\alpha_3, \ldots \}
\endaligned
\tag 1.3
$$
coincides as a multiset (i.e., up to the order) with the list of
row frequencies in {\rm(1.1)}. Moreover, $\widetilde{\beta}_k=\beta_k$
for all $k=1,2,\ldots$.
\endproclaim

\smallskip\noindent{\bf Example.}
In the particular case of $\alpha_1=1$, $\alpha_k=0$ for $k\ge2$,
$\beta_k=0$ for $k\ge1$, the row frequencies form a geometric
sequence
$$
\widetilde{\alpha} = \{(1-t),(1-t)t,(1-t)t^2,\ldots\}.
\tag 1.4
$$
In this case the conjecture follows from a result obtained by
A.~Borodin \cite{Bor}.
\medskip

Let $M$ be a probability distribution on the compact group $U(k)$
of unipotent matrices. Following the general terminology of the
theory of graded graphs and locally semisimple algebras (cf.\
\cite{LSS}),
we say that the measure $M$ is {\it central} if the
probabilities
$$
M_\lambda \overset{\text{\rm def}}\to{=}
M\{u\in U(k):\quad u(n)=v\}
\tag 1.5
$$
depend not on the matrix $v\in U_n(k)$ itself but on its Jordan
form partition $\lambda$ only. We reduce the main problem of
describing the characteristic measures to the more general
problem of classifying all central probability measures on
$U(k)$. As usual, every such measure admits a unique presentation
as a mixture of ergodic (indecomposable) central measures. The
latter can be described in terms of the (dual) Hall--Littlewood
symmetric polynomials $Q_\lambda(x;t)$ (see \cite{Mac}).

\proclaim{Conjecture}
Given a real number $b$, we write
$b^{(t)}=\{(1-t)b,t(1-t)b,t^2(1-t)b,\ldots\}$ for the {\it
associated geometric sequence, and we set
$\beta^{(t)}=\{\beta_1^{(t)},\beta_2^{(t)},\ldots\}$ for a
sequence $\beta$. Then the formula
$$
M_\lambda = (1-t)^{-n}t^{-n(\lambda)}Q_\lambda(\alpha;\beta^t;t)
\tag 1.6
$$
correctly defines the cylinder probabilities {\rm(1.5)} of a central
measure $M=M^{(\alpha;\beta)}$, where $(\alpha;\beta)$ is a point
of the Thoma simplex. The measures $M^{(\alpha;\beta)}$ are ergodic
and mutually singular, and every ergodic central measure $M$
coincides with one of the measures $M^{(\alpha;\beta)}$.
\endproclaim

This conjecture is equivalent to the conjecture of the second
author stated in \cite{GHL, equation (7.3.3)}
(see also \cite{Ker}).

\head Unipotent characters of the groups $GL_n(k)$ \endhead

\subhead Unipotent representations of the group $GL_n(k)$
\endsubhead
An irreducible representation of the group $G_n=GL_n(k)$ (and its
character) is called {\it unipotent} if it is a part of some flag
representation (that is, a representation induced by a trivial
representation of a parabolic subgroup). Equivalently, one can
define unipotent representations as those containing an invariant
vector for the Borel subgroup $B_n\subset G_n$ of upper
triangular matrices.

The unipotent representations $\pi_\lambda$ are labelled by the
integer partitions $\lambda\vdash n$ of $n$ (we denote the set
of such partitions by $\Bbb Y_n$). It is well known that the dimension of
$\pi_\lambda$ is given by the $q$-hook formula
$$
\dim \pi_\lambda = \frac
{(q-1)(q^2-1)\ldots(q^n-1)}{\prod_{b\in\lambda} (q^{h(b)}-1)}
\;q^{n(\lambda)} .
$$
The dimension of the subspace of $B_n$-invariant vectors does not
depend on $q$ and is provided by the ordinary hook formula
$$
f_\lambda = \frac{n!}{\prod_{b\in\lambda} h(b)}.
$$

\smallskip\noindent
{\bf Example.}
The so-called {\it Steinberg representation} $\pi_{(1^n)}$ is the
unipotent representation of $G_n$ corresponding to the one-column
partition $\lambda=(1^n)$. Its dimension is
$$
\dim \pi_{(1^n)} = q^{n(n-1)/2},
$$
and there is only one $B_n$-invariant vector.
\medskip

Let $\chi^\lambda(g)=$ Tr$\,\pi_\lambda(g)$ denote the character
of an irreducible unipotent representation $\pi_\lambda$ of the
group $G_n$. There is a simple explicit formula for the values
$\chi^\lambda_\rho(q)$ of such characters at unipotent elements
$g\in G_n$ of type $\rho$:
$$
\chi^\lambda_\rho(q) = \widetilde{K}_{\lambda\rho}(q) \equiv
q^{n(\rho)}\, K_{\lambda\rho}(q^{-1}).
$$
Recall that $K_{\lambda\rho}(t)$ is the generalized Kostka
matrix, see \cite{Mac, Section III.6}. In particular, all the
values $\chi^\lambda_\rho(q)$ are nonnegative.

Alternatively, one can define the character matrix
$\{\chi^\lambda_\rho(q)\}$ as the transition matrix between the
bases of Schur functions $s_\lambda(x)$ and that of modified Hall
polynomials $\widetilde{P}_\rho(x;q)\equiv
q^{-n(\rho)}P_\rho(x;1/q)$:
$$
s_\lambda(x) = \sum_{\rho\vdash n} \chi^\lambda_\rho(q)\,
\widetilde{P}_\rho(x;q).
$$
Note the similarity of this formula with the classical Frobenius
formula providing the character matrix of the symmetric group.

\subhead Representations induced from parabolic subgroups
\endsubhead
Given a composition $\mu=(\mu_1,\mu_2,\,\ldots,\mu_m)$, we denote
by $P_\mu$ the corresponding parabolic subgroup generated by the
block diagonal matrices in
$G_{\mu_1}\times G_{\mu_2}\times\ldots\times G_{\mu_m}$ and by
the group $B_n$. Let $\sigma_\mu$ be the representation of $G_n$
induced by the identity representation of the subgroup
$P_\mu\subset G_n$. We shall refer to $\sigma_\mu$ as  {\it
induced} representations. Note that the induced representation
$\sigma_\mu$ depends only on the partition corresponding to the
composition $\mu$ (that is, does not depend on the order of
parts). By definition, $\sigma_\mu$ is the representation in the
space of flags $H_0=\{0\}\subset H_1\subset\ldots\subset
H_m=k^n$ with the dimension vector
$\big(\dim(H_j)-\dim(H_{j-1})\big)_{j=1}^m=\mu$. The dimension of
the induced representation $\sigma_\mu$ equals the Gaussian
coefficient:
$$
\dim \sigma_\mu = \frac{[n]!}{\prod_{j=1}^m [\mu_j]!},
$$
where $[n]!=\prod_{j=1}^n[j]$ and $[j]=(q^j-1)/(q-1)$.

The important fact is that the multiplicities $K_{\lambda,\mu}$
of unipotent representations in the decomposition of a flag
representation,
$$
\sigma_\mu = \sum_\lambda K_{\lambda,\mu} \pi_\lambda,
$$
do not depend on $q$, and coincide with the {\it Kostka numbers}
arising in decompositions of induced characters of the symmetric
group $\frak S_n$.

Let $\psi^\mu(g)=$ Tr$\,\sigma_\mu(g)$ denote the character
of the induced representation $\sigma_\mu$. By definition,
$\psi^\mu(g)$ is simply the number of fixed points in the space
of the representation $\sigma_\mu$ (i.e., of flags of a specified
dimension vector $\mu$). The value of the induced character
$\psi^\mu$ at a unipotent element $g\in G_n$ of type $\rho$ can
be written in the form
$$
\psi^\mu(g) = \psi^\mu_\rho(q) \equiv \sum_{\lambda\vdash n}
K_{\lambda,\mu}\, \widetilde{K}_{\lambda,\rho}(q).
$$
This is a polynomial in $q$.

\subhead The characters of unipotent factor representations
\endsubhead
We refer to \cite{VK} for the description of characters of
general factor representations of the group $GLB(k)$ (more
precisely, of the group algebra $A$ of Bruhat--Schwartz functions
on this group).

We say that such a character $\chi$ is {\it unipotent}
\footnote{In \cite{VK} we have used the term {\it principal
series} instead.} if its restriction $\chi\big|_{A_n}$ to every
subalgebra $A_n\cong\Bbb C[G_n]$ is a linear combination of only
unipotent characters of the group $G_n$. According to
\cite{VK, Section 6}, the unipotent characters of
factor representations $\chi^{(\alpha;\beta)}$ are labelled by
the points $(\alpha;\beta)\in \Delta$ of the Thoma simplex
$\Delta$. We shall use two explicit decompositions of these
characters: one in terms of irreducible unipotent characters
$\chi^\lambda$,
$$
\chi^{(\alpha;\beta)}(a_g) = \sum_{\lambda\vdash n}
\chi^\lambda(g)\, s_\lambda(\alpha;\beta),
$$
and one in terms of the characters $\psi^\mu$ of induced representations,
$$
\chi^{(\alpha;\beta)}(a_g) = \sum_{\nu\vdash n}
\psi^\nu(g)\, m_\nu(\alpha;\beta).
$$

\subhead Induced characters at primary elements \endsubhead
The goal of this section is to describe the value of an induced
character $\psi^\mu$, $\mu\vdash n$, at a primary element $g\in
G_n$ in terms of its values at unipotent elements. Recall (see
\cite{Mac, IV.3}) that the primary
conjugacy class corresponding to
an irreducible polynomial $f\in k[t]$ of degree $d=d(f)$ is
characterized by an integer partition $\rho\vdash m$, $m=n/d$. If
$d=1$ and $f=t-1$, the class is said to be unipotent.

\proclaim{Lemma}
The value $\psi^\mu(g)$ of the character $\psi^\mu$ at a primary element $g\in G_n$
with the characteristic polynomial $f^n(t)$, $d(f)=1$, does not
depend on the choice of $f$.
\endproclaim

\demo{Proof}
Let $g\in G_n$ be a primary element with the characteristic
polynomial $f^n(t)$ where $f(t)=t-a$. Then $g=aE+n$, where $E\in
G_n$ is the identity matrix and $n\in G_n$ is a nilpotent
element. Let $h\in G_n$ be another primary element with
the same partition $\mu\vdash n$ but a different polynomial
$t-b$, $b\in k$. We can assume that $h=bE+n$. Since $g-h=(a-b)E$
is a scalar operator, $g$ and $h$ have identical invariant
subspaces $V\subset k^n$. Therefore, $g$ and $h$ have the same
number of invariant flags of every dimension type $\mu$ and
$\psi^\mu(g)=\psi^\mu(h)$.
\qed\enddemo

We are now in a position to consider primary elements with an
irreducible polynomial $f$ of degree $d>1$. Let $K\cong k^d$ be
the extension field of $k$ of degree $d$. The polynomial $f$ has
$d$ distinct roots in $K$, and we denote by $a\in K$ one of the
roots. The operator of multiplication by $a$ in $K$ can also be
treated as an operator in $k^d$. It is well known that this is a
primary element of $GL_d(k)$ with the characteristic polynomial
$f$ and the integer partition $\mu=(1)$. Every semisimple regular
element is conjugate to some multiplication operator in the
extension field (see \cite{Mac, VI.3, Example 4}).

More generally, let $g\in GL_m(K)$ be a primary element with a
partition $\mu\vdash m$ and a linear polynomial $t-a$ where $a\in
K$ generates the extension $K\supset k$. Denote by
$\widetilde{g}\in GL_n(k)$ the corresponding operator in $k^n$,
$n=md$. Then $\widetilde{g}$ is also primary with the same
integer partition $\mu$. The corresponding polynomial $f$ of
degree $d$ can be identified as the irreducible polynomial over
$k$ with the root $a$. Every primary element in $GL_n(k)$ can be
obtained in this way from a primary element with a polynomial of
degree $d=1$ over an appropriate extension field.

One can easily describe the invariant subspaces of
$\widetilde{g}$ in $k^n$ in terms of those of $g$ in $K^m$.

\proclaim{Lemma}
Let $H\subset K^m$ be a $g$-invariant subspace of dimension $s$
over $K$. Then it is also a $\widetilde{g}$-invariant subspace
over $k$ of dimension $sd$. Moreover, every
$\widetilde{g}$-invariant subspace in $k^n$ is in fact a subspace
over the extension field $K$. In particular, its $k$-dimension is
divisible by $d$.
\endproclaim

\proclaim{Corollary}
Let $g\in GL_n(k)$, $n=md$, be a primary element with partition
$\rho$. Then the value $\psi^\nu(g)$ of an induced character
$\psi^\nu$ only depends on the degree $d=d(f)$ of the irreducible
polynomial $f$ corresponding to $g$. More precisely,
$$
\psi^\nu(g) = \cases
\psi^\mu_\rho(q^d) &
\text{ if } \nu = d\mu \text{ for some } \mu\vdash m \\
0 & \text{ otherwise}.
\endcases
$$
\endproclaim

Recall that $d\mu$ denotes the partition of $n=md$ obtained from
$\mu\vdash m$ by multiplying all its parts by $d$. The
polynomial $\psi^\mu_\rho(q)$ was defined at the end of Section
2.2 as the value of a character on a unipotent element.

\subhead An extension formula for unipotent characters of the
group $GLB(k)$ \endsubhead
Let $\chi^{\alpha;\beta}$ be the unipotent character of the group
$GLB(k)$ corresponding to Thoma parameters $(\alpha;\beta)$.
According to the formulas of Section 2.3, the value of
$\chi^{\alpha;\beta}$ at a unipotent element $g\in GL_n(k)$ with
partition $\rho$ can be written in the form
$$
\chi^{\alpha;\beta}(g) = r_\rho(\alpha;\beta;t),
$$
where $t=1/q$ and the symmetric function $r_\rho$ is defined by
the formula
$$
r_\rho(\cdot) = t^{-n(\rho)} \sum_{\lambda\vdash n}
K_{\lambda,\rho}(t)\; s_\lambda(\cdot).
$$

In this section we find an explicit expression for the value of a
unipotent character $\chi^{\alpha;\beta}$ at a general element
$g\in GL_n(k)$ in terms of the symmetric functions $r_\rho$.

\proclaim{Theorem}
Assume that the conjugacy class of $g\in GL_n(k)$ corresponds to
a family $\varphi:\,\Cal F\to\Bbb Y$ such that $\sum_f
\deg(f)\,|\varphi(f)|=n$ (see \cite{Mac, IV.2}). Then
$$
\chi^{\alpha;\beta}(g) = \prod_{f\in\Cal F}
r_{\varphi(f)}(\alpha^{d(f)};-(-\beta)^{d(f)};t^{d(f)}).
$$
\endproclaim

\demo{Proof}
Recall (see \cite{Mac, I.2}) that the monomial symmetric function
$m_\lambda(x_1,x_2,\ldots)$ is defined as the sum of the monomials
$x^\gamma=x_1^{\gamma_1}\,x_2^{\gamma_2}\ldots$ where $\gamma$
runs over all distinct finite permutations of the coordinates of
$\lambda$. Let $\pi_d$ denote the endomorphism of the symmetric
function algebra $\Lambda$ defined by $\pi_d(p_m)=p_{md}$, for
all $m=1,2,\ldots$. This is the familiar {\it plethysm}
endomorphism, see \cite{Mac, I.8}. Then it follows immediately from the
definitions that $\pi_d(m_\lambda)=m_{d\lambda}$.

Given a symmetric function $u\in\Lambda$, set $U=\pi_d(u)$. By
definition, we have $U(x)=u(x^d)$. We
shall use the following ``supersymmetric'' version of this
formula:
$$
U(\alpha;\beta) = u(\alpha^d;-(-\beta)^d).
$$
In fact, it suffices to check the formula for the power sum
symmetric functions $u=p_m$. One can easily check that in this case both
sides are equal to $\sum \alpha_i^{md}-\sum(-\beta_i)^{md}$.

Let us now prove the theorem in the special case of a primary
element $g$ with a polynomial of degree $d$. By the Lemma of section
2.4,
$$
\aligned
\chi^{\alpha;\beta}(g)
&= \sum_{\nu\vdash n} \psi^\nu(g)\, m_\nu(\alpha;\beta) =\\
&= \sum_{\mu\vdash m} \psi^{d\mu}(g)\, m_{d\mu}(\alpha;\beta) =\\
&= \sum_{\mu\vdash m} \psi^\mu_\rho(q^d)\,
m_\mu(\alpha^d;-(-\beta)^d) =\\
&= r_\rho(\alpha^d;-(-\beta)^d;t^d).
\endaligned
$$

The theorem now follows from the Multiplication Theorem, see
\cite{VK, Section 7, Theorem 6}.
\qed\enddemo

\subhead Unipotent characters at unipotent classes \endsubhead
Let $\chi=\chi^{(\alpha;\beta)}$ be the unipotent character of
the group $GLB(k)$ corresponding to a point $(\alpha;\beta)$ of
the Thoma simplex. Let $\rho$ be a Jordan form partition of a
matrix $u\in GL_n(k)$. In order to describe the character value
$\chi(u)$ we introduce appropriate symmetric functions
$$
r_\rho(\alpha;\beta) = t^{-n(\rho)} \sum_{\lambda\vdash n}
K_{\lambda,\rho}(t)\; s_\lambda(\alpha;\beta),
$$
where $s_\lambda$ is the Schur function. See \cite{VK} and \cite{Mac, Section
III.6} for the definition of the generalized Kostka matrices
$K_{\lambda,\rho}(t)$.

\proclaim{Theorem}
The value $\chi^{(\alpha;\beta)}(u)$ of a
unipotent character $\chi^{(\alpha;\beta)}$
at a unipotent
matrix $u$ with Jordan form partition $\rho$ can be written as
$$
\chi^{(\alpha;\beta)}(u) = r_\rho(\alpha;\beta).
$$
These characters of $GLB(k)$ are determined by their
values at unipotent classes.
\endproclaim

\newpage

\topmatter

\title
IV. An outline of construction of 
factor representations of the
group $GLB(\Bbb F_q)$ \endtitle
\endtopmatter

\subhead {1. Introduction}
\endsubhead
It is well known that indecomposable semifinite characters of the
group $S_{\infty}$,  which were described in \cite{FA},
correspond to type $II_\infty$ representations of the group 
algebra ${\Bbb C}(S_{\infty})$ (and, consequently, the group itself).
Recall that in this case, in general, a character takes infinite
values on the group. Since the algebra $A$ (the Bruhat--Schwartz algebra
for the group $GLB(\Bbb F_q)$, see \cite{VK}) is a two-sided invariant
dense subalgebra of the true group algebra $L^1(GLB,\mu)$
of the group  $GLB(\Bbb F_q)$ with the Haar measure $\mu$, 
it follows by a general
theorem (see, e.g., \cite{HR}) that 
representations
of the algebra $A$ can be extended to the group $GLB$ and, consequently, to the group
 $GL(\infty)$, which is a subgroup of $GLB$. But, in general, this extension
 can have another type; more exactly, if a representation
 of $A$ is of type $II$, then its extension can be either
of the same type or of type $III$. It is this possibility 
that occurs in our case.
 However, even without referring to the general theory, which,
 as far as we know, gives no
 answer to the question how one can find the type of the extension, 
 we can nevertheless
 give a direct description of some representations of the group
 corresponding to type $II_\infty$ 
 representations of the algebra $A$.

 We will restrict ourselves to one example, namely, the important
 class of representations which we call {\it principal series
 Grassmanian representations}. The described 
 construction is also of independent interest.

\subhead 2. The infinite Grassmanian, Schubert cells, and their symbols
\endsubhead
In this section we generalize the classical methodology of
constructing the cell partition of the Grassmanian of 
finite-dimensional subspaces into Schubert cells (see, e.g., \cite{Miln,
Stan}) to the infinite-dimensional case. More precisely, we construct
the partition of the Grassmanian of infinite-dimensional
subspaces into Schubert cells.

\define\Gr{\operatorname{Gr}}

 Consider the Grassmanian $\Gr(V)$ of all subspaces of a 
 vector space $V$ over the finite field $\Bbb F_q$, where
 $V=\lim V_n$,  $V_1 \subset V_2\subset \dots$,  $\dim {V_n}=n$,
 and endow it with the topology of a $GLB(q)$-homogeneous space. In other words, a
 neighborhood of a given subspace is the set of subspaces
 that can be obtained from it by applying
 elements of the group $GLB$ lying in a given
 neighborhood of the identity element. It is convenient to choose
 the fundamental base
 of neighborhoods of a given subspace $t \in \Gr(V)$
 that consists of the sets \{$U_{t,n}$\}, where ${U_{t,n}}$
 is the set of subspaces $s \in \Gr(V)$ that have the same intersection
 with the coordinate subspace $V_n$ as $t$ and
 the same dimensions of the intersections with  $V_m$, $m \geq n$, as $t$,
 i.e., the set of subspaces such that
 $s \cap V_n=t\cap V_n$ and
 $\dim (s \cap V_m) = \dim(t \cap V_m)$, $m \geq n$; we will say that
 such a subspace $s$ is $n$-equivalent to $t$.
 This topology  on the
 Grassmanian is stronger than the totally disconnected topology. 
 The flag of subspaces  $V_1 \subset V_2\subset \dots$, $\dim
 {V_n}_{F_q}=n$, $n=1,2, \dots$, will be called the {\it principal flag}.

 Every subspace $E \in \Gr(V)$, $E \subset V$, is uniquely determined by
 the infinite sequence $E \cap V_i$, $i=1,2, \ldots$, of its (finite-dimensional) intersections
 with the subspaces of the principal flag.
 Consider the dimensions  $\dim(E \cap V_i)=d_i$; the sequence
 $\epsilon_i=d_i-d_{i-1}$, $i=1,2, \ldots$, of the differences of
 neighboring dimensions, which are equal to $0$ or $1$,
 is the {\it Schubert symbol} of the subspace $E$
 (it is convenient for us to slightly change the traditional terminology).
 The subspaces with the same Schubert symbol generate a {\it Schubert cell},
 which, as it is easy to see, is compact 
 and is a
 $B$-orbit,  where $B$ is the compact Borel subgroup
 of the group $GLB$, i.e., the group of upper triangular matrices.
 Thus the Grassmanian $\Gr(V)$  is a fibration over the space of
 $(0,1)$-sequences, with Schubert cells as fibers.
 Every fiber (cell) can be identified with an affine
 space of some dimension over the field $\Bbb F_q$.
 The dimension of a cell $e$ is equal to $\dim(e)={\sum_{i=1}^\infty
 i\epsilon_i}$. Cells of infinite dimension consist of
 infinite-dimensional subspaces. We are most interested in cells 
 that consist of subspaces of infinite dimension and infinite
 codimension.

 Every Schubert cell contains a unique coordinate subspace, which
 can be used as a distinguished point of this cell.
Given a $(0,1)$-sequence $\epsilon$, denote by $S_\epsilon$ the
Schubert cell with symbol $\epsilon$, and by $t_\epsilon$ the
distinguished (coordinate) subspace that belongs to this cell.
Denote by $\Gr_0(V)$ the set of coordinate subspaces and
identify it with the space $(0,1)^ \infty \equiv \Sigma$ of all
$(0,1)$-sequences by the formula $$\epsilon \mapsto t_\epsilon. $$ 
Thus the set of  $(0,1)$-sequences parameterizes
simultaneously Schubert symbols, coordinate subspaces,  and
Schubert cells.

Since the group $B$ is compact and all Schubert cells are
transitive $B$-spaces, on every Schubert cell 
there is a unique $B$-invariant probability
measure, namely, the image of the Haar measure on
$B$. Denote it by $M_\epsilon$.

{\it The orbit of a Schubert cell} $S$ with respect to the action of the
countable group $GL(\infty)$ is the countable union of Schubert
cells over all symbols that differ from the symbol of
$S$ in finitely many coordinates and have
the same number of $1$'s in any sufficiently large initial part of
the symbol. Now the orbit is a $GLB$-invariant set, because the group $GLB$
is generated by the subgroups $B$ and $GL(\infty)$. The
partition of the space of symbols into the classes of symbols 
belonging to the same orbit is exactly the partition $\tau$ of the
space $\Sigma = (0,1)^\infty $ into the orbits of the group
$S_\infty$; it is a tame (hyperfinite) partition. Two subspaces $E,E'
\in \Gr(V)$ that belong to Schubert cells of the same orbit
(in other words, whose symbols eventually coincide) 
will be called {\it congruent}.

On the orbit of every cell $S_\epsilon=S$ we have the $\sigma$-finite
$GLB$-invariant measure $\mu_S$  normalized by the
condition that the measure of the given cell $S_\epsilon$ (or of the
given $B$-orbit) is equal to $1$. We will call such measures {\it
elementary invariant measures}. The measures corresponding to
different cells of the same orbit differ by a positive factor.
More exactly, the following proposition holds.

\proclaim{Proposition} The elementary invariant measure of a
compact set is finite, while the measure of any cylinder
in $\Gr(V)$ is infinite. Let $S$ and $S'$ be two
cells whose symbols $\epsilon$ and $\epsilon'$ coincide
starting from the $n$th coordinate, and let
$\mu_{S}$ and $\mu_{S'}$ be the corresponding elementary measures.
Then
$$   \mu_{S}(C)/\mu_{S'}(C)=q^{\dim(S \cap V_n)-\dim(S' \cap V_n)}=
   q^{\sum_{i=1}^\infty i(\epsilon_i - \epsilon'_i)},  $$
where $C$ is an arbitrary compact set. Thus the ratio
of the measures of cells is equal to the ratio of the exponentials of
their dimensions.
\endproclaim

\medskip\noindent{\bf Remark.}
The last sum is finite because the symbols coincide starting
from the $n$th coordinate. Denote the expression introduced in the
proposition above by $c(S,S')$; it is defined on pairs of 
Schubert cells from the same orbit and, consequently, on pairs of
subspaces $E,E'$ belonging to such cells: $c(E,E')=c(S,S')$, $E \in
S$, $E' \in S'$.
\medskip

The elementary $\sigma$-finite invariant measures are parameterized
by the classes of Schubert cells (or by the elements of the
partition $\tau$).
The function $c(\cdot,\cdot)$
will be called the {\it fundamental cocycle} on the Grassmanian $\Gr(V)$.
It is also possible to regard elementary invariant measures
mentioned above as the images of the Haar measure of the group
 $GLB$ on orbits of cells; however, we cannot directly
project the Haar measure onto an orbit (regarded as a homogenous space),
because the stationary (parabolic) subgroup is of infinite Haar
measure.

Using this transitive action of the group $GLB$ with the invariant
$\sigma$-finite measure $\mu_S$, we can, as usual,
define a unitary representation of the group $GLB$ in the space
$L^2(\mu_S)$. But we will study more complicated {\it
factor representations} of $GLB$.

Let us define new measures on the Grassmanian. First consider a
measure $\nu$ on the space of symbols $\Sigma$. Taking the direct
integral
   $$ \int M_\epsilon d\nu(\epsilon) \equiv M^{\nu} $$
over this measure, we obtain a $B$-invariant measure $M^{\nu}$ on
the Grassmanian. Let us find conditions under which this measure is
quasi-invariant
with respect to the group $GLB$ and central with respect to the
algebra $A$ (see \cite{VK}).

\proclaim{Proposition} The measure $M^\nu$
is $GLB$-quasi-invariant if and only if
the measure $\nu$ is quasi-invariant 
with respect to the natural action of the group $S_\infty$ on the
space $\Sigma$. The action of $GLB$ is ergodic if and only if
the measure $\nu$ is ergodic with respect to
the action of $S_\infty$, or, equivalently, if the tail
partition $\tau$ is ergodic.
\endproclaim

Having the {\it congruence relation} on the
Grassmanian $\Gr(V)$, we can define the principal groupoid 
corresponding to this relation (see \cite{Ren}); denote it by $R$.
 The measure $M^{\nu}$ determines a measure $\mu^{\nu}$  on the groupoid
regarded  as a subset of $\Gr(V) \times \Gr(V)$ (see\ \cite{Ren}).
  Let $H=L^2(R, \mu^{\nu})$ be the Hilbert space of complex
  square-integrable (with respect to this measure) functions on the groupoid.
  According to general constructions (see \cite{LSS}),
unitary representations of the groups $GL(\infty)$ and $GLB$ (and, consequently,
$*$-representations of the algebra $A$) are defined in this Hilbert space
with the help of a left ($ \pi_l$) and a 
right ($ \pi_r$) actions, respectively. 

First  consider the representation of the algebra $A$ in the space $H$.
Let $a_g \in A$; then $\pi_l(a_g)$ is an operator that is defined
by a kernel and acts as the operator of left multiplication
by a function $f_g$ on $R$.

Now let us introduce an ergodic measure on the set of Schubert symbols
that is invariant with respect to the partition $\tau$;
according to the remark above, such a measure must be
$S_\infty$-invariant,  and, consequently, by the well-known de
Finetti theorem, it is a Bernoulli measure with parameters $(\alpha,
1-\alpha)$, $0\le \alpha \le 1$; denote it by $\mu^\alpha$.

Consider the direct integral 
of the $B$-invariant measures on the Schubert cells
with respect to this measure
$\mu^\alpha$;
by definition, this is a Bernoulli measure on the Grassmanian;
 we denote it by $\mu_\alpha$. This measure is $B$-invariant
 and $GLB$-quasi-invariant by definition.

The orbits of the group $GLB$ coincide with the 
unions of Schubert cells mentioned above.
Let us define the Radon--Nikodym (R-N) cocycle of the measure $\mu_\alpha$
with respect to the group $GLB$. It is more convenient to
regard the cocycle as a function of pairs of points $E$ and $E'$ on the
same orbit. Let $E$ and $E'$ be two subspaces
belonging to the same Schubert cell; this means that there exists
an element $g\in GLB$ that sends $E$ to
$E'$. Then, by definition, the value of the cocycle 
at the pair $E,E'$ is equal to the value of
the density $d\mu_\alpha(gE)/d\mu_\alpha$.  
The symbols $\{\epsilon_i(E)=\epsilon_i\}$ and
$\{\epsilon_i(E')=\epsilon'_i\}$
coincide starting from some
coordinate with index  $n(E,E')=n$, and the numbers of $1$'s among
the coordinates of these symbols with indices less than $n$ are the same and
equal to $k(E,E')=k$. Since the action of the Borel subgroup
preserves the measure and since this action is transitive on the
cells, the value of the cocycle depends on the cell only, i.e., the cocycle
is a function of symbols and can be calculated on the
coordinate subspaces.

\proclaim{Proposition} Let $e$ and $e'$ be the intersections of the
subspaces $E$ and $E'$, respectively, with $V_n$, and let their dimensions
be equal to $k$. Let $\sigma(e)$ and $\sigma(e')$ be their Schubert
cells in the Grassmanian $\Gr_k(n)$. Then the R-N-cocycle of the measure
$\mu_\alpha$ with respect to the action of the group $GLB$ is equal
to
$$c(E,E')=q^{\dim(\sigma(e))-\dim(\sigma(e'))}=
   q^{\sum_{i=1}^\infty i(\epsilon_i - \epsilon'_i)},  $$
where $\dim$ is the dimension of a finite-dimensional Schubert cell
regarded as a manifold over the field $\Bbb F_q$.
\endproclaim

\proclaim {Corollary} For any $n$, the restriction of the cocycle $c$ to the
subgroup $GL_n(\infty)$ is cohomologous to one. At the
same time, the cocycle  $c$ is not cohomologous to zero on the whole group
$GLB$, hence there is no (finite or $\sigma$-finite) $GLB$-invariant
measure equivalent to $\mu_\alpha$.
\endproclaim

\proclaim{Proposition} The action of the group $GLB$ 
on the space $\Gr(V)$ with the
quasi-invariant measure $\mu_\alpha$ is
ergodic.
\endproclaim

\demo {Proof} Since the action of the group $B$ is transitive, it
suffices to make sure that the action of a subgroup of
$GLB(\infty)$ on the space $(\Gr_0(V), \mu^\alpha)$ of coordinate
subspaces is ergodic. Consider the subgroup $S_\infty$ of 
$GLB(q)$; its action is ergodic on $(\Gr_0(V),\mu^\alpha)$, because
it coincides with the action of
$S_\infty$ on the space $\Sigma$ of all
$(0,1)$-sequences with a Bernoulli measure, which is ergodic by
the well-known zero or one law.
\enddemo

\medskip\noindent{\bf Remark.}
For every $n$ there exists a $GL_n(\infty)$-invariant measure that is
equivalent to $\mu_\alpha$.
\medskip

\subhead 3. Several remarks on the realization of principal series
factor representations of the groups
$GLB(q)$ and $GL(\infty,q)$
\endsubhead
A brief sketch of the construction of factor representations
  is as follows.

The construction is based on the same idea that was used in our papers  \cite{VK81, VK82},
namely, the idea of trajectory (or orbit, or groupoid) model of
representations, which goes back to von Neumann. But this
construction uses the orbit equivalence relation, obtained from
the orbit partition, instead of the action of groups. We
start with introducing a $GLB(q)$-quasi-invariant and
$B$-invariant measure on the Grassmanian $\Gr(V)$, and then, using
the orbit equivalence relation on $\Gr(V)$ and this measure, we
construct the Hilbert space $L^2$ in which the required
factor representation of the group $GLB(q)$ is realized in a natural way.

  Consider the $GLB(q)$-orbit equivalence relation $\tau$
 on  $(\Gr(V),\mu_\alpha)$  and construct the Hilbert space
 $L^2(\tau)$ (see \cite{VK81}). In this space we have a left and a right
 unitary representations of the group $GLB(q)$ and, consequently, of the
 algebra $A$. It turns out that these representations of
 $GLB(q)$ are of type $III$, while their restrictions to the algebra
 $A$ are of type $II_\infty$; the characters coincide with the
 characters described in the previous draft.
 It seems that the phenomenon of changing the type of factors
 when extending a representation as was mentioned above had not been
 considered in detail in representation theory.

\subhead 4. Realization of 
factor representations related to the Grassmanian
\endsubhead
The goal of this section is to describe the factor representation of the
algebra $A$ that corresponds to a principal series character
$\chi_\alpha$ with two frequencies $\alpha=(\alpha_1,\alpha_2)$. 
These representations are constructed from a measure on the Grassmanian.
The
construction below is slightly different from the construction in 
Sec.~3 of this draft.

We start from a graph interpretation of the {\it mixed grassmanian} $\Gr$
of all possible $k$-linear subspaces of the space $V$. Given a
subspace $E\in \Gr$, denote by $E_n=E\bigcap V_n$ its intersection
with the $n$th space $V_n$ of the principal flag (see above). For
$n=1,2,\ldots$, set $d_n(E)=\dim E_n-\dim E_{n-1}$. The sequence
$d_n=d_n(E)$ consists of $0$'s and $1$'s and determines a path in the
Pascal triangle. We draw this graph in the SE quadrant.
The value
$d_n=0$ corresponds to a horizontal edge (in this case 
$E_n=E_{n-1}$), and the value
$d_n=1$ corresponds to a vertical edge (the dimension of $E_n/E_{n-1}$ is
equal to $1$). If the subspace $E_n$ is generated by the subspace $E_{n-1}$
and the coordinate vector $e_n$, then the extension is
called a {\it coordinate extension}. Of course, besides the coordinate
extension, there are $q^m-1$ other extensions, where $m$ is the
codimension of $E_{n-1}$ in $V_{n-1}$ (as well as the codimension of $E_n$ in
$V_n$). With each such extension we associate the vertical edge
between the corresponding vertices of the Pascal triangle.

The {\it Pascal $q$-triangle} is the branching graph $P(q)$ whose
vertices and edges are the same as in the ordinary Pascal
triangle $P=P_1$, but the vertical edges $(i-1,j)\to(i,j)$ from the
$j$th column are assigned the multiplicity $q^j$. The dimensions in
the graph $P(q)$ are nothing else than the Gaussian binomial
coefficients, and the paths of length $n$ from the
initial vertex parameterize our subspaces $E_n\subset V_n$.

The Grassmanian $\Gr$ is identified with the space of infinite paths
in the (multi)graph $P(q)$ (recall that a path is a
sequence of edges). The group $GLB(q)$ acts in a natural way on
$V$ and, consequently, on $\Gr(V)$. Orbits of the Borel
group $B$ (i.e., Schubert cells, see above) correspond to {\it plaits}, i.e.,
sets of paths with the same vertices. Every plait contains
a distinguished path corresponding to the sequence of 
coordinate subspaces. The set of distinguished paths determines a
subgraph of $P(q)$ which is the ordinary Pascal triangle $P$, so we
have the inclusion $P(1)\subset P(q)$.

Let us say that a subspace $E$ is {\it almost coordinate} 
if almost all extensions
$E_{n-1}\subset E_n$  are coordinate extensions; let $\Gr_n$ be the
subset of subspaces $E\subset V$ that are coordinate starting
from some index $n$. Let $\Gr_\infty\bigcup \Gr_n$ be the set of all
almost coordinate subspaces $E\subset V$.

We endow the space $\Gr_n$ with a finite measure $M_n$ which  is defined as
the average of the Bernoulli measure with parameters
$(\alpha_1,\alpha_2)$ on the set of coordinate subspaces with respect
to the action of the finite group $G_n$ permuting the first $n$ edges
of paths in the graph $P(q)$. The direct description of this
measure is as follows. Let each vertical edge have the weight
$\alpha_1$ and each horizontal edge have the weight $\alpha_2$. Then the
$M_n$-measure of the cylinder of paths with common initial part $t$ is equal
to the product of the weights of the edges that belong to $t$. The
full measure of the space $\Gr_n$ is equal to
$$
M_n(\Gr_n) = \sum_{m=0}^n {n \choose m}_q \alpha_1^m \alpha_2^{n-m},
$$
where ${n\choose m}_q$ are the Gaussian binomial coefficients. When
we extend the Grassmanian $\Gr_n$ to $\Gr_{n+1}$, the measure
$M_{n+1}$ extends the measure $M_n$. Thus we have
a $\sigma$-finite infinite measure $M$
 on the space
$\Gr_\infty=\Gr(V)$. This measure
coincides with one of the measures constructed in the
previous section; it possesses invariant and quasi-invariant
properties with respect to the groups $B$, $GL(\infty)$, $GLB$ and,
consequently, determines the corresponding factor representations.

\enddocument